\documentclass[10pt]{article}
\usepackage{amssymb,amsmath,color}
\usepackage [hypertex]{hyperref}

\allowdisplaybreaks
\begin{document}

\def\pd#1#2{\frac{\partial#1}{\partial#2}}
\def\dfrac{\displaystyle\frac}
\let\oldsection\section
\renewcommand\section{\setcounter{equation}{0}\oldsection}
\renewcommand\thesection{\arabic{section}}
\renewcommand\theequation{\thesection.\arabic{equation}}

\def\Xint#1{\mathchoice
  {\XXint\displaystyle\textstyle{#1}}%
  {\XXint\textstyle\scriptstyle{#1}}%
  {\XXint\scriptstyle\scriptscriptstyle{#1}}%
  {\XXint\scriptscriptstyle\scriptscriptstyle{#1}}%
  \!\int}
\def\XXint#1#2#3{{\setbox0=\hbox{$#1{#2#3}{\int}$}
  \vcenter{\hbox{$#2#3$}}\kern-.5\wd0}}
\def\ddashint{\Xint=}
\def\dashint{\Xint-}

\newcommand{\cblue}{\color{blue}}
\newcommand{\cred}{\color{red}}
\newcommand{\be}{\begin{equation}\label}
\newcommand{\ee}{\end{equation}}
\newcommand{\bea}{\begin{eqnarray}\label}
\newcommand{\eea}{\end{eqnarray}}
\newcommand{\nn}{\nonumber}
\newcommand{\intO}{\int_\Omega}
\newcommand{\Om}{\Omega}
\newcommand{\cd}{\cdot}
\newcommand{\pa}{\partial}
\newcommand{\ep}{\varepsilon}
\newcommand{\uep}{u_{\eta}}
\newcommand{\vep}{v_{\eta}}
\newcommand{\Sep}{S_{\eta}}
\newcommand{\Lp}{L^p(\Om)}
\newcommand{\Lq}{L^q(\Om)}
\newcommand{\abs}{\\[2mm]}
\newcommand{\eps}{\varepsilon}
\newtheorem{thm}{\indent Theorem}
\newtheorem{lem}{\indent Lemma}[section]
\newtheorem{proposition}{\indent Proposition}[section]
\newtheorem{dnt}{\indent Definition}[section]
\newtheorem{remark}{\indent Remark}[section]
\newtheorem{cor}{\indent Corollary}[section]
\allowdisplaybreaks

\title{Multiple solutions for a class of Kirchhoff equation with singular nonlinearity }
\author{Zupei Shen$^{\rm\,a}$,~ Zhiqing Han$^{\rm a,\,}${\thanks{Corresponding author.
E-mail: hanzhiq@dlut.edu.cn (Z.-Q. Han), pershen@126.com (Zupei.Shen).Tel(Fax): +86 411
84707268.}} }
\date{$^{\rm a}${\small School of Mathematical Sciences, Dalian University of Technology, Dalian 116024, PR China}\\}
\maketitle
\date{}
{\bf Abstract}. In this article, we investigate the existence and multiplicity of solutions of  Kirchhoff equation
\begin{equation*}
 \left\{
\begin{aligned}
-(1+b \int_{\mathbb{R}^3}|\nabla u|^2)\Delta u=  k(x)\frac{|u|^2 u}{|x|} +\lambda h(x)u,~~~~x\in\mathbb{R}^3\\
u(x)\rightarrow 0  ~~~~~~~~~~~~~~~~~~~~~~~~~as~~ |x|\rightarrow\infty
\end{aligned}
\right.
\end{equation*}
where the potential $k(x)$ allows sign changing. Making use of Nehari manifold method and concentration-compactness principle,
we obtain the existence and multiplicity of solutions for this equation. Our main results can be viewed as  partial extensions of the results of \cite{C-K-W, Chen, Liu-Sun}.

\section{Introduction}
\mbox{}\indent
The system
\begin{equation}
 \left\{
\begin{aligned}
-(a+b \int_{\mathbb{R}^3}|\nabla u|^2)\Delta u +V(x)u= f(x,u)~~~~x\in\mathbb{R}^3\\
u(x)\rightarrow 0  ~~~~~~~~~~~~~~~~~~~~~~~~~as~~ |x|\rightarrow\infty \label{11}
\end{aligned}
\right.
\end{equation}
is related to the stationary analogue of the equation
\begin{equation}
 \rho\frac{\partial^2 u}{\partial t^2}-(\frac{p_0}{h}+\frac{E}{2L}\int^L_0|\frac{\partial u}{\partial t}|^2 dx)\frac{\partial^2 u}{\partial x^2}=0, \label{12}
\end{equation}
which was presented by  Kirchhoff in 1883. Kirchhoff's model takes into account the changes in length of the string produced by
transverse vibrations. The parameters in (\ref{12}) have following physical meanings:
$L$ is the length of the string, $h$ is the area of the cross-section, $E$ is the Young modulus of the material, $\rho$  is the
mass density and $P_0$ is the initial tension.
After J. L. Lions's work \cite{lions1}, which introduced an abstract functional analysis framework to the following equation
\begin{equation}
u_{tt}-(a+b\int_\Omega|\nabla u|^2)\triangle u=f(x,u) \label{13},
\end{equation}
equation (\ref{11}) received much attention. See \cite{Ancona,Arosio,Wang} and the references therein.
A typical way to deal with  equation (\ref{11}) is to use the mountain pass theorem \cite{A-R}. For this purpose, one usually assumes that $f(x,u)$ is subcritical, superlinear at the origin and satisfies the Ambrosetti-Rabinowitz condition (AR in short): e.g. see \cite{Chen-li}
$$\exists \mu>4 ~~such ~that~ 0<\mu F(x,u)\leq f(x,u)u ~for ~all ~u\in \mathbb{R}.  $$
Using Nehari manifold method, He and Zou \cite{He}  proved the existence of positive ground state solution of (\ref{11}) with the nonlinearity satisfying the Ambrosetti-Rabinowitz condition. The typical case is $f(u)\sim |u|^{p-2}u $ with $4<p<6$.
Wu \cite{Wu} obtained the existence of nontrivial solutions to a class of Kirchhoff equation. He assumed that the nonlinearity $f(x,u)$ is 4-superlinear at infinity and satisfies
$$4F(x,u)\leq f(x,u) ~\text{for ~all }~u\in \mathbb{R}.$$
In order to get compactness, he considered the problem in a weight subspace
$$E\triangleq\{u\in H^1|\int_{\mathbb{R}^3}V(x)|u|^2 dx<\infty\}$$
such that $ E\hookrightarrow L^p$ is compact.
Li and Ye \cite{Li-Ye} partially extended  the results of He and Zou to $3<p<6$ by  monotonicity trick and a global compactness lemma.
There are also some works to deal with Kirchhoff equation with indefinite nonlinearity. Recently, Chen, Kuo and Wu  \cite{C-K-W} investigated the multiplicity of positive solutions for the problem which involving sign-changing weight functions
 \begin{equation}
 \left\{
\begin{aligned}
-(a+b\int_{\mathbb{R}^3}|\nabla u|^2)\Delta u&=  k(x){|u|^{p-2} u} +\lambda h(x){|u|^{q-2} u} ,~~~~x\in\Omega\\
u(x)&= 0  ~~~~~~~~~~~~~~~~~~~~~~~~~~~~~~~~~~~~~~~~~~~~~x\in \partial \Omega
\end{aligned}
\right. \label{100}
\end{equation}
where $\Omega$ is a smooth bounded domain in $\mathbb{R}^3$ with $1<q<2<p<6$. The authors showed that existence and multiplicity of results strongly depend on the size of $p$  with respect to 4. Part of the results is the following:
 If $ p=4$, then the problem has (at least) one solution for $b$ large and two positive solutions for $b$ and $\lambda$ small.
Chen \cite{Chen}  proved that equation
\begin{equation}
 \left\{
\begin{aligned}
-(1+b \int_{\mathbb{R}^3}|\nabla u|^2)\Delta u=  k(x){|u|^{p-2} u} +\lambda h(x)u,~~~~x\in\mathbb{R}^3\\
u(x)\rightarrow 0  ~~~~~~~~~~~~~~~~~~~~~~~~~as~~ |x|\rightarrow\infty
\end{aligned}
\right. \label{15}
\end{equation} exists  multiple positive solutions,
where  $k(x)$  allows sign changing with $p\in(4,6)$.
As for  singular nonlinearity, Liu and Sun  \cite{Liu-Sun}  considered the existence of positive solutions for the following problem with singular and superlinear terms
\begin{equation}
 \left\{
\begin{aligned}
-(a+b \int_{\mathbb{R}^3}|\nabla u|^2)\Delta u &= h(x)u^{-r} +\lambda k(x)\frac{|u|^{p-2} u}{|x|^s} ~~~~x\in\Omega\\
u(x)&=0  ~~~~~~~~~~~~~~~~~~~~~~~~~~~~~~~~~~~x\in \partial \Omega \label{10}
\end{aligned}
\right.
\end{equation}
where $0\leq s<1,~~~~ 4<p<6-2s,~~~0<r<1 ~\text{and }~ k(x)\geq0.$  They obtained two positive solutions by Nehari manifold.
However, very little is known for existence of nontrivial solutions of (\ref{11}) if $f(x,u)$ is  singular and indefinite.
Motivated by \cite{Liu-Sun,C-K-W, Chen}, in the present paper, we consider the case where $f(u,x)$ is a combination of a singular 4-linear term and a linear term. More precisely, we study the following system with the form
\begin{equation}
 \left\{
\begin{aligned}
-(1+b \int_{\mathbb{R}^3}|\nabla u|^2)\Delta u=  k(x)\frac{|u|^2 u}{|x|} +\lambda h(x)u,~~~~x\in\mathbb{R}^3\\
u(x)\rightarrow 0  ~~~~~~~~~~~~~~~~~~~~~~~~~as~~ |x|\rightarrow\infty
\end{aligned}
\right.  \label{ks}
\end{equation}
where $b>0, ~h(x)>0 $ and $k(x)$ is indefinite. In order to state our main results, we assume the following hypotheses (H):
\begin{itemize}
\item[(H$_{h}$)] $h\in L^\frac{3}{2}(\mathbb{R}^3)$, $ h(x)\geq 0$ for any $x\in \mathbb{R}^3 $;
\item[(H$_{k_1}$)] $k(x)\in C(\mathbb{R}^3)$ and $k(x)$ changes sign in $\mathbb{R}^3$;
\item[(H$_{k_2}$)] $\lim_{|x|\rightarrow\infty}k(x)=k_\infty < 0, k(0)=0$.
\end{itemize}
As far as  we know, no one considered this case before. Under  hypothesis (H$_{h}$), there exists a sequence of eigenvalues $\lambda_n$ of
$$-\Delta u+u=\lambda h(x)u ~~~~~~~~~\text{in}~~ H^1(\mathbb{R}^3) $$
with $0 < \lambda_1 < \lambda_2 \leq\cdots$ and each eigenvalue being of finite multiplicity. The associated normalized eigenfunctions are denoted
by $e_1, e_2 \cdots $ with $\|e_i\|=1$. Moreover, $e_1 > 0$ in $ \mathbb{R}^3 $.

We are now ready to state our results:
\begin{thm}{\label{thm1}}
Assume that hypotheses (H) hold. Then for $0 <\lambda <\lambda_1$, problem (\ref{ks}) has at least one solution in $D^{1,2}(\mathbb{R}^3)$.
\end{thm}
\begin{thm}{\label{thm2}}
 Assume that hypotheses (H) hold and  $\int_{\mathbb{R}^3}\frac{k(x)}{|x|}{e_1}^4dx-b(\int_{\mathbb{R}^3}|\nabla e_1 |^2dx)^2 < 0$. Then there exists $\delta > 0$ such that problem (\ref{ks}) has  at least two solutions whenever $\lambda_1<\lambda <\lambda_1+\delta$.
\end{thm}
\begin{flushleft}
  \emph{\textbf{Remark 1}}. Comparing with  problem (\ref{100}), we mainly consider the problem in the whole space ${\mathbb{R}^3}$ with $q=2,p=4 $. In this sense, our main results can be viewed as  partial extensions of the results of \cite{C-K-W}.\\
\end{flushleft}
\emph{\textbf{Remark 2}}.  To the best of our knowledge, for the semilinear elliptic equations with indefinite nonlinearity, a similar condition like  $\int_{\mathbb{R}^5}k(x){e_1}^q dx< 0$ is needed (e.g.see \cite{Alama, C-V, Brown}).
In \cite{Chen}, the authors proved similar results for equation (\ref{15}) as $4<p<6$, and the condition  $\int_{\mathbb{R}^3}\frac{k(x)}{|x|}{e_1}^4-b(\int_{\mathbb{R}^3}|\nabla e_1 |^2dx)^2 < 0$ was not needed.  However, he needed another  condition:
\begin{itemize}
\item[(A2)]: $ |\Omega^0 |=0 ~\text{where} ~\Omega^0 = \{x\in \mathbb{R}^3 :k(x)=0\}$.
\end{itemize}
By using the same argument in this paper, it is much easier to get the same result for equation (\ref{15}) when $p=4$ . In this sense, our main results can be viewed as a partial extension of the result of \cite{Chen}.\\
\emph{\textbf{Remark 3}}. Comparing with problem (\ref{10}), we mainly consider the case that $p=4, s=1$. In this sense, our main results can be views as a partial extension of the results of \cite{Liu-Sun}.\\
 \emph{\textbf{Remark 4}}. For system (1.1), when the nonlinearity is subcritical, as far as we know, no one consider the "zero mass" case, that is $V(x)=0$.\\
To prove Theorems 1.1 and  1.2, we use the Nehari manifold method borrowing from Brown and Zhang \cite{Brown}. In \cite{Brown}, the authors considered a semilinear boundary value problem on a bounded domain. Together with a concentration-compactness principle, Chabrowski and Costa \cite{Chabrowski-Costa, Chabrowski2} generalized the result to unbounded region and singular nonlinearity respectively.
Inspired by the  papers of Brown -Zhang \cite{Brown} and  Chabrowski-Costa \cite{Chabrowski-Costa, Chabrowski2}, we extend the results to the  Kirchhoff equation in $\mathbb{R}^3$.

\section{Preliminaries}
To go further, let us give some notions and some  known results.
\begin{itemize}
\item[*] $D^{1,2}(\mathbb{R}^3)$ is the completion of $C_0^\infty(\mathbb{R}^3)$ with respect to the norm $\| u\|^2 = \int_{\mathbb{R}^3}|\nabla u |^2dx$.
\item[*] $\|\cdot\|$ denotes the norm of $D^{1,2}(\mathbb{R}^3)$.
\item[*] $\rightarrow$ denotes the strong convergence.
\item[*] $\rightharpoonup$ denotes the weak convergence.
\item[*] $C, C_i$ and $c$ denote various positive constants.
\end{itemize}
\subsection{$C^1$~functional}
For $u\in D^{1,2}(\mathbb{R}^3)$, weak solutions to (\ref{ks}) correspond to critical points of the energy functional

$$J(u)=\frac{1}{2}\int_{\mathbb{R}^3}|\nabla u|^2-\lambda h(x)u^2dx + \frac{b}{4}(\int_{\mathbb{R}^3}|\nabla u|^2dx)^2-\frac{1}{4}\int_{\mathbb{R}^3}\frac{k(x)| u|^4}{|x|}dx.$$
By the Caffarelli-Kohn-Nirenberg inequality \cite{C-K-N}
\begin{equation*}
C(\int_{\mathbb{R}^3}|x|^{-qb}|u|^q)^{\frac{p}{q}}\leq \int_{\mathbb{R}^3}|x|^{-pa}|\nabla u|^pdx.
\end{equation*}
where $ u\in C_0^\infty(\mathbb{R}^3)$, $ 1<p<N $ ,  $0\leq a\leq b\leq a+1\leq\frac{N}{P}$,   $ q:= \frac{NP}{N+p(b-a)-p}$.
Let $u\in D^{1,2}(\mathbb{R}^3)$, by approximation, it is easy to see that there exists a constant C such that
\begin{equation}
C(\int_{\mathbb{R}^3}|x|^{-1}|u|^4)^{\frac{1}{2}}\leq \int_{\mathbb{R}^3}|\nabla u|^2dx.\ \label{14}
\end{equation}
(Let $a=0, ~b=\frac{1}{4}, ~N=3 ,~p=2$, we have $q:=p_*=4$.)\\
By (\ref{14}), it is no difficult  to show that the functional $J$ is of class $C^1$ (See Lemma 2.2). Moreover,$$J'(u)v=\int_{\mathbb{R}^3}\nabla u~\nabla v-\lambda h(x)uvdx +\int_{\mathbb{R}^3}|\nabla u|^2dx \int_{\mathbb{R}^3}\nabla u~\nabla vdx  -\int_{\mathbb{R}^3}\frac{k(x)}{|x|}|u|^2uvdx$$
 for any $v\in D^{1,2}(\mathbb{R}^3)$.

In order to use the critical point theory, we need to prove  the energy functional~$J(u)$ is of a class of $C^1$ in $D^{1,2}(\mathbb{R}^3)$ .
\begin{lem}
If $u_n \rightarrow u$~in $D^{1,2}(\mathbb{R}^3)$,~ there exists a subsequence, still denoted by $u_n$ and $g\in D^{1,2}(\mathbb{R}^3)$ such that $u_n \rightarrow u$ almost everywhere on $\mathbb{R}^3$ and $$|u_n|\leq g,|u|\leq g.$$
\end{lem}
\emph{Proof.} Going if necessary to a subsequence, we can assume that $u_n \rightarrow u$ $a.e.$ on $\mathbb{R}^3$. There exists a subsequence, still denoted by $u_n$ such that
$$\|u_{j+1}-u_j\|\leq 2^{-j},~~~~~\forall~~j\geq 1.$$
Let us define
$$g(x)=|u_1(x)|+\sum_1^\infty|u_{j+1}(x)-u_j(x)|.$$
It is clear that $|u_n|\leq g,|u|\leq g$  a.e. on $\mathbb{R}^3$ and $g\in D^{1,2}(\mathbb{R}^3).$\\
\emph{\textbf{Remark}}. In order to use the Lebesgue convergence theorem and  Caffarelli-Kohn-Nirenberg inequality, in the proof of Lemma 2.2 we require $g\in D^{1,2}(\mathbb{R}^3)$. Usually, we only require $g\in L^p(\mathbb{R}^3)$ for $1\leq p <\infty$; see  Lemma A.1 of Appendix A in \cite{Willem}. In this sense, the lemma seems to be  new.
\begin{lem}
 $J$ is of class $C^1$ in $D^{1,2}(\mathbb{R}^3).$
\end{lem}
\emph{Proof.}
Let  $\varphi(x)=\int_{\mathbb{R}^3}\frac{u^4}{|x|}dx.$
We only need to prove $\varphi(u)$ is of $C^1$ class.\\
First, we prove to \emph{Existence of the Gateaux derivative of $\varphi$ at $u$ }.\\
Let $u$, $v\in  D^{1,2}(\mathbb{R}^3)$ and $t\in (0,1)$. Since
$$(u+tv)^4=u^4+C_4^1u(tv)^3+C_4^2u^2(tv)^2+C_4^3u^3(tv)^1+(tv)^4.$$
We have
\begin{eqnarray*}
\frac{|((u+tv)^4-u^4)|}{|x|t}&\leq &\frac{4t^2|u||v|^3+6t|u|^2|v|^2+4|u|^3|v|+t^3|v|^4}{|x|}\\
   &\leq &6\frac{|u||v|^3+|u|^2|v|^2+|u|^3|v|+|v|^4}{|x|}.
\end{eqnarray*}
By H\"{o}lder inequality and Caffarelli-Kohn-Nirenberg inequality (\ref{14}), we have
\begin{eqnarray*}
\int_{\mathbb{R}^3}\frac{|u||v|^3}{|x|}dx&=&\int_{\mathbb{R}^3}\frac{|u|}{|x|^{\frac{1}{4}}}\frac{|v|^3}{|x|^{\frac{3}{4}}}dx\\
&\leq& (\int_{\mathbb{R}^3}(\frac{|v|^3}{|x|^\frac{3}{4}})^{\frac{4}{3}})^\frac{3}{4}(\int_{\mathbb{R}^3}(\frac{|u|}{|x|^\frac{1}{4}})^4)^\frac{1}{4}\\
&\leq& C(\int_{\mathbb{R}^3}|\nabla v|^2)dx)^\frac{3}{2}(\int_{\mathbb{R}^3}|\nabla u|^2)dx)^\frac{1}{2}\\
&\leq& C_1.
\end{eqnarray*}
Similarly,
$$\int_{\mathbb{R}^3}\frac{|u|^2|v|^3}{|x|}dx\leq C_2, \int_{\mathbb{R}^3}\frac{|u|^3|v|}{|x|}dx\leq C_3, \int_{\mathbb{R}^3}\frac{|v|^4}{|x|}dx\leq C_4. $$
Therefore, we have
$$\frac{|((u+tv)^4-u^4)|}{|x|t}\leq \eta(v)\in L^1(\mathbb{R}^3)$$ where $\eta(v)=6\frac{|u||v|^3+|u|^2|v|^2+|u|^3|v|+|v|^4}{|x|}$.
It follows then from the  Lebesgue convergence theorem that
$$\langle\varphi'(u),v\rangle=4\int_{\mathbb{R}^3}\frac{|u|^2uv}{|x|}dx.$$
Next, we prove the \emph{ Gateaux derivation is continuous}.\\
Assume  that  $u_n \rightarrow u$~in $D^{1,2}(\mathbb{R}^3)$. By H\"{o}lder inequality and Caffarelli-Kohn-Nirenberg inequality (\ref{14}), we have
\begin{eqnarray*}
|\langle\varphi'(u_n)-\varphi'(u),v\rangle|&=& \int_{\mathbb{R}^3}\frac{(|u_n|^2u_n-|u|^2u)v}{|x|}dx\\
&=&\int_{\mathbb{R}^3}\frac{|v|}{|x|^\frac{1}{4}}\frac{(|u_n|^2u_n-|u|^2u)v}{|x|^\frac{3}{4}}dx\\
&\leq& C\left(\int_{\mathbb{R}^3}\left(\frac{(|u_n|^2u_n-|u|^2u)v}{|x|^\frac{3}{4}} \right)^\frac{4}{3}dx\right)^\frac{3}{4}\|v\|.
\end{eqnarray*}
Lemma 2.1 implies that there exists $g\in D^{1,2}(\mathbb{R}^3)$ such that
$|u_n|\leq g,|u|\leq g$. By  Caffarelli-Kohn-Nirenberg inequality (\ref{14}), we have
$$\left(\frac{(|u_n|^2u_n-|u|^2u)}{|x|^\frac{3}{4}} \right)^\frac{4}{3}\leq \frac{2g^4}{|x|}\in L^1(\mathbb{R}^3).$$
 According to  Lebesgue convergence theorem, we get
 $$\|\varphi'(u_n)-\varphi'(u)\|\leq C\left(\int_{\mathbb{R}^3}\left(\frac{(|u_n|^2u_n-|u|^2u)}{|x|^\frac{3}{4}} \right)^\frac{4}{3}\right)^\frac{3}{4}\rightarrow 0 ~~ \text {as}~~ n\rightarrow \infty.$$
\subsection{Nehari manifold}
For $u\in D^{1,2}(\mathbb{R}^3)$, weak solutions to (\ref{ks}) correspond to critical points of the energy functional

$$J(u)=\frac{1}{2}\int_{\mathbb{R}^3}|\nabla u|^2-\lambda h(x)u^2dx + \frac{b}{4}(\int_{\mathbb{R}^3}|\nabla u|^2dx)^2-\frac{1}{4}\int_{\mathbb{R}^3}\frac{k(x)| u|^4}{|x|}dx.$$

 Since the functional $J$ is not bounded from below on $D^{1,2}(\mathbb{R}^3)$, a good candidate for an appropriate subset to study $J$ is the so-called Nehari manifold
\begin{eqnarray*}
S&=&\{{u\in D^{1,2}(\mathbb{R}^3)}~| ~J'(u)u=0\}\\
&=& \{{u\in D^{1,2}(\mathbb{R}^3)}~ |\int_{\mathbb{R}^3}|\nabla u|^2-\lambda h(x)u^2dx=\int_{\mathbb{R}^3}\frac{k(x)| u|^4}{|x|}dx-(\int_{\mathbb{R}^3}|\nabla u|^2dx)^2\}.
\end{eqnarray*}
It is useful to understand $S$ in  term of the stationary points of the fibering mappings, i.e. $$\varphi_u(t)= J(tu)=\frac{t^2}{2}\int_{\mathbb{R}^3}|\nabla u|^2-\lambda h(x)u^2dx-\frac{t^4}{4}\int_{\mathbb{R}^3}\frac{k(x)| u|^4}{|x|}dx+\frac{bt^4}{4}(\int_{\mathbb{R}^3}|\nabla u|^2dx )^2.$$
We now follow some ideas from the paper \cite{Brown}.
\begin{lem}
 Let $ u\in D^{1,2}(\mathbb{R}^3)-\{0\}$ and $t > 0$. Then $tu\in S$ if and only if $\varphi'_u(t)=0$.
\end{lem}
Thus the points in $S$ correspond to the  stationary points of the fiber map $\varphi_u(t)$ and so it is natural to divide $S$ into three parts $S^+, S^-$ and $S^0$ corresponding to local minima, local maxima and points of inflexion of the fibering maps. We have$$\varphi''_u(t)=\int_{\mathbb{R}^3}|\nabla u|^2-\lambda h(x)u^2dx-3t^2(\int_{\mathbb{R}^3}\frac{k(x)| u|^4}{|x|}dx -b(\int_{\mathbb{R}^3}|\nabla u|^2dx )^2)$$\\
and
$$\varphi''_u(1)=\int_{\mathbb{R}^3}|\nabla u|^2-\lambda h(x)u^2dx-3(\int_{\mathbb{R}^3}\frac{k(x)| u|^4}{|x|}dx -b(\int_{\mathbb{R}^3}|\nabla u|^2dx )^2)dx.$$
Hence if we define

$$S^+=\{u\in S:\int_{\mathbb{R}^3}|\nabla u|^2-\lambda h(x)u^2dx-3(\int_{\mathbb{R}^3}\frac{k(x)| u|^4}{|x|}dx -b(\int_{\mathbb{R}^3}|\nabla u|^2dx )^2)dx >0\},$$
$$S^-=\{u\in S:\int_{\mathbb{R}^3}|\nabla u|^2-\lambda h(x)u^2dx-3(\int_{\mathbb{R}^3}\frac{k(x)| u|^4}{|x|}dx -b(\int_{\mathbb{R}^3}|\nabla u|^2dx )^2)dx < 0\},$$
$$S^0=\{u\in S:\int_{\mathbb{R}^3}|\nabla u|^2-\lambda h(x)u^2dx-3(\int_{\mathbb{R}^3}\frac{k(x)| u|^4}{|x|}dx -b(\int_{\mathbb{R}^3}|\nabla u|^2dx )^2)dx=0\},$$
we have
\begin{lem}{\label{lem2}}
 Let $u\in S.$ Then
\begin{eqnarray*}
S^+& =&\{u\in S:\int_{\mathbb{R}^3}|\nabla u|^2-\lambda h(x)u^2dx<0\}\\
   &=&\{u\in S:\int_{\mathbb{R}^3}\frac{k(x)| u|^4}{|x|}dx -b(\int_{\mathbb{R}^3}|\nabla u|^2dx )^2 < 0\},
\end{eqnarray*}
\begin{eqnarray*}
S^-& =&\{u\in S:\int_{\mathbb{R}^3}|\nabla u|^2-\lambda h(x)u^2dx>0\}\\
   &=&\{u\in S:\int_{\mathbb{R}^3}\frac{k(x)| u|^4}{|x|}dx -b(\int_{\mathbb{R}^3}|\nabla u|^2dx )^2 > 0\},
\end{eqnarray*}
\begin{eqnarray*}
S^0& =&\{u\in S:\int_{\mathbb{R}^3}|\nabla u|^2-\lambda h(x)u^2dx =0\}\\
   &=&\{u\in S:\int_{\mathbb{R}^3}\frac{k(x)| u|^4}{|x|}dx -b(\int_{\mathbb{R}^3}|\nabla u|^2dx )^2 = 0\}.
\end{eqnarray*}
\end{lem}
Since $$\varphi'_u(t)=t\int_{\mathbb{R}^3}|\nabla u|^2-\lambda h(x)u^2dx-t^3(\int_{\mathbb{R}^3}\frac{k(x)| u|^4}{|x|}dx-b(\int_{\mathbb{R}^3}|\nabla u|^2dx )^2),$$ $\varphi_u$ has exactly one  turning point at $t(u)=\{\frac{\int_{\mathbb{R}^3}|\nabla u|^2-\lambda h(x)u^2dx}{\int_{\mathbb{R}^3}\frac{k(x)| u|^4}{|x|}dx-b(\int_{\mathbb{R}^3}|\nabla u|^2dx )^2}\}^\frac{1}{2}$ if and only if $\int_{\mathbb{R}^3}|\nabla u|^2-\lambda h(x)u^2dx$ and $\int_{\mathbb{R}^3}\frac{k(x)| u|^4}{|x|}dx-b(\int_{\mathbb{R}^3}|\nabla u|^2dx )^2$ have the same sign.\\
As in \cite{Brown}, we let
$$L^+=\{u\in H^1(\mathbb{R}^3) :\|u\|=1,\int_{\mathbb{R}^3}|\nabla u|^2-\lambda h(x)u^2dx > 0\},$$%-3
$$L^-=\{u\in H^1(\mathbb{R}^3) :\|u\|=1,\int_{\mathbb{R}^3}|\nabla u|^2-\lambda h(x)u^2dx < 0\},$$
$$L^0=\{u\in H^1(\mathbb{R}^3) :\|u\|=1,\int_{\mathbb{R}^3}|\nabla u|^2-\lambda h(x)u^2dx = 0\},$$
and
$$B^+=\{u\in H^1(\mathbb{R}^3) :\|u\|=1,\int_{\mathbb{R}^3}\frac{k(x)| u|^4}{|x|}dx-b(\int_{\mathbb{R}^3}|\nabla u|^2dx )^2> 0\},$$
$$B^-=\{u\in H^1(\mathbb{R}^3) :\|u\|=1,\int_{\mathbb{R}^3}\frac{k(x)| u|^4}{|x|}dx-b(\int_{\mathbb{R}^3}|\nabla u|^2dx )^2 < 0\},$$
$$B^0=\{u\in H^1(\mathbb{R}^3) :\|u\|=1,\int_{\mathbb{R}^3}\frac{k(x)| u|^4}{|x|}dx-b(\int_{\mathbb{R}^3}|\nabla u|^2dx )^2 = 0\}.$$

\begin{lem}{\label{lem3}}
(i). A multiple of $u$ lies in $S^-$ if and only if $\frac{u}{\|u\|}$ lies in $ L^+\cap B^+.$\\
(ii). A multiple of $u$ lies in $S^+$ if and only if $\frac{u}{\|u\|}$ lies in $L^-\cap B^-.$\\
(iii). For  $u\in L^+\cap B^- $ or $u\in L^-\cap B^+ $, no multiple of $u$ lies in $S$.
\end{lem}

\begin{thm}{\label{thm3}}
Suppose that $u_0$ is a local minimizer  for $J$ on $S$ and $u_0\notin S^0$, then $J'(u_0)=0.$
\end{thm}
\subsection{concentration-compactness principle}
In order to overcome the loss of  compactness we make use of a simple version of concentration-compactness principle which is  considered in \cite[\text{Proposition 1.3}]{Chabrowski1}, see also \cite[\text{Concentration-Compactness Principle}]{Chabrowski2}.
\begin{thm}
Let $(u_m)$ be a sequence  in $D^{1,2}(\mathbb{R}^3)$ such that
\begin{eqnarray*}
u_m(x)\rightarrow u(x) ~~~a.e. ~~in~~ \mathbb{R}^3,\\
u_m(x)\rightharpoonup u(x)~~~in ~~D^{1,2}(\mathbb{R}^3),\\
|\nabla(u_m-u)|^2 \rightharpoonup \mu ~~in ~~M(\mathbb{R}^3),\\
|x|^{-\frac{1}{4}}|u_m-u|^4\rightharpoonup \nu~~in ~~M(\mathbb{R}^3),\\
\end{eqnarray*}
where $M(\mathbb{R}^3)$ denotes the space of bounded measures in $\mathbb{R}^3$. Define the quantities
\begin{equation*}
\alpha_\infty=\lim_{R\rightarrow \infty}\limsup_{n\rightarrow \infty}\int _{|x|>R}\frac{|u_n|^4}{|x|}dx
\end{equation*}
\begin{equation*}
\beta_\infty=\lim_{R\rightarrow \infty}\limsup_{n\rightarrow \infty}\int _{|x|>R}|\nabla u_n|^2 dx.
\end{equation*}
\begin{equation*}
\nu_0=\lim_{R\rightarrow \infty}\limsup_{n\rightarrow \infty}\int _{|x|\leq \frac{1}{R}}\frac{|u_n|^4}{|x|}dx
\end{equation*}
\begin{equation*}
\mu_0=\lim_{R\rightarrow \infty}\limsup_{n\rightarrow \infty}\int _{|x|\leq \frac{1}{R}}|\nabla u_n|^2 dx.
\end{equation*}
Then we have
\begin{equation}
\limsup_{n\rightarrow \infty}\int_{\mathbb{R}^3}|\nabla u_n|^2dx =\int_{\mathbb{R}^3}|\nabla u|^2 dx +\beta_\infty +\mu_0.  \label{23'}
\end{equation}
\begin{equation}
\limsup_{n\rightarrow \infty}\int_{\mathbb{R}^3}\frac{| u_n|^4}{|x|}dx=\int_{\mathbb{R}^3}\frac{| u|^4}{|x|}dx+\alpha_\infty + \nu_0.  \label{24'}
\end{equation}
\end{thm}
\section{The case when $0 <\lambda <\lambda_1$}
Suppose that  $ 0 <\lambda <\lambda_1$. It is easy to see that there exists $ \theta>0$ such that
\begin{equation}
\int_{\mathbb{R}^3}|\nabla u|^2-\lambda h(x)u^2dx \geq \theta \|u\|^2. \label{31}
\end{equation}
Thus $S^+$ is empty and $S^0=\{0\}$.\\
To prove the theorem, we need the following lemma.

\begin{lem}{\label{lem3.1}}Suppose $ 0 <\lambda <\lambda_1$. Then\\
(i).  $\inf_{u\in S^-}J(u)>0$;\\
(ii). There exists $u \in S^- \setminus \{0\}$, such that $J(u) =\inf_{v\in S^-}J(v)$.
\end{lem}
\emph{Proof.}
(i). By Lemma 2.4, we have$$J(u)=\frac{1}{4}\int_{\mathbb{R}^3}k(x)|u|^4-l(x){\phi_u(x)}u^2dx > 0 ~~~~ when ~~~u\in S^-.$$ So $J$ is bounded below by $0$ on $S^-$.
We  show that $\inf_{u\in S^-}J(u)>0$. Suppose $u\in S^-$. Then $v=\frac{u}{\|u\|}\in L^+\bigcap B^+$ and $u=t(v)v$ where $t(v)=[\frac{\int_{\mathbb{R}^3}|\nabla v|^2-\lambda h(x)v^2dx}{\int_{\mathbb{R}^3}k(x)\frac{|v|^4}{|x|}dx -b(\int_{\mathbb{R}^3}|\nabla u|^2dx)^2 }]^\frac{1}{2}$.
In addition,
\begin{eqnarray*}
J(u)=J(t(v)v)&=&\frac{1}{4}(t(v))^2\int_{\mathbb{R}^3}|\nabla v|^2+|v|^2-\lambda h(x)v^2dx\\
&=&\frac{1}{4}\frac{{(\int_{\mathbb{R}^3}|\nabla v|^2-\lambda h(x)v^2dx)}^2}{\int_{\mathbb{R}^3}k(x)\frac{|v|^4}{|x|}dx-b(\int_{\mathbb{R}^3}|\nabla u|^2dx)^2 }\\
&\geq& \frac{\theta^2}{4\int_{\mathbb{R}^3}k(x)\frac{|v|^4}{|x|}dx-b(\int_{\mathbb{R}^3}|\nabla u|^2dx)^2)}~~~~~~~(by~~ (\ref{31}))\\
&\geq& \frac{\theta^2}{4\int_{\mathbb{R}^3}k(x)\frac{|v|^4}{|x|}dx}.~~~~~~~~~~~~~~~~~~~~~~~~~~~~~~~~~~~~~~~~~~~~~~~~~~~~~~~~~~~~~~~~~~~~~~~~~~~~~~~~~~~~~~~~~~~~~~~~(a)
\end{eqnarray*}
We now focus on the term $\int_{\mathbb{R}^3}k(x)\frac{|v|^4}{|x|}dx$.   By Caffarelli-Kohn-Nirenberg inequality (\ref{14}),
we have
\begin{eqnarray*}
\int_{\mathbb{R}^3}k(x)\frac{|v|^4}{|x|}dx&\leq & C|k|_{L^\infty}\int_{\mathbb{R}^3}\frac{|v|^4}{|x|}dx\\
&\leq&C|k|_{L^\infty}\|v\|^2\\
&=& C|k|_{L^\infty}.~~~~~~~~~~~~~~~~~~~~~~~~~~~~~~~~~~~~~~~~~~~~~~~~~~~~~~~~~~~~~~~~~~~~~~~~~~~~~~~~~~~~~~~~~~~~~~~~~~~~~~~~~~~~~~~~(b)
\end{eqnarray*}
Combining ($a$) with ($b$), we have $$J(u)\geq\frac{\theta^2}{C|k|_{L^\infty}}>0.$$ Hence $$\inf_{u\in S^-}J(u)>0.$$
(ii). We show that there exists a minimizer on $S^-$. Let $\{u_n\}\subset$ $S^-$ be a minimizing sequence, i.e, $\lim_{n\rightarrow\infty}J(u_n) =\inf_{u\in S^-}J(u)$. By ({\ref{31}}), we have
\begin{eqnarray*}
J(u_n)=\frac{1}{4}\int_{\mathbb{R}^3}|\nabla u_n|^2-\lambda h(x)u_n^2dx \geq\frac{1}{4}\theta\|u_n\|^2.
\end{eqnarray*}
So $\{u_n\}$ is bounded in $D^{1,2}(\mathbb{R}^3)$. Passing to a subsequence if necessary, we obtain that $\{u_n\}\rightharpoonup u $ in $D^{1,2}(\mathbb{R}^3)$. First, we claim that $u\neq 0$. Since $\{u_n\}\subset$ $S$, we have
\begin{eqnarray*}
\int_{\mathbb{R}^3}|\nabla u_n|^2-\lambda h(x)u_n^2dx =\int_{\mathbb{R}^3}k(x)\frac{|u_n|^4}{|x|}dx-b(\int_{\mathbb{R}^3}|\nabla u_n|^2dx)^2.
\end{eqnarray*}
Using (\ref{23'}) and (\ref{24'}), we deduce
\begin{equation}
\int_{\mathbb{R}^3}|\nabla u|^2-\lambda h(x)|u|^2dx +\beta_\infty +\mu_0 \leq \int_{\mathbb{R}^3}k(x)\frac{|u|^4}{|x|}dx-b(\int_{\mathbb{R}^3}|\nabla u|^2dx)^2+k(0)\nu_0+k(\infty)\nu_\infty. \label{32}
\end{equation}

Suppose $u=0$.  By (\ref{32}), we have $$\beta_\infty=\mu_0 =0.$$
Then $u_n \rightarrow 0$ in $D^{1,2}(\mathbb{R}^3)$, a contradiction to $\inf_{u\in S^-}J(u)>0$.\\
We now claim that  $\beta_\infty=\mu_0 =0$. Otherwise, we deduce from (\ref{32}) that
$$0<\int_{\mathbb{R}^3}|\nabla u|^2-\lambda h(x)|u|^2dx < \int_{\mathbb{R}^3}k(x)\frac{|v|^4}{|x|}dx-b(\int_{\mathbb{R}^3}|\nabla u|^2dx^2)dx.$$
There exists $0 < s < 1$ such that
$$\int_{\mathbb{R}^3}|\nabla su|^2-\lambda h(x)|su|^2dx = \int_{\mathbb{R}^3}\frac{k(x)}{|x|}(su)^4dx-b(\int_{\mathbb{R}^3}|\nabla su|^2dx^2).$$
This implies that $su$ belongs to $S^-$. On other hand, since
\begin{align*}
0 <\int_{\mathbb{R}^3}|\nabla u|^2-\lambda h(x)|u|^2dx &\leq \liminf_{n\rightarrow\infty}\int_{\mathbb{R}^3}|\nabla u_n|^2-\lambda h(x)|u_n|^2dx\\
&=4\inf_{w\in S^-}J(w)\\
&\leq \int_{\mathbb{R}^3}|\nabla su|^2-\lambda h(x)|su|^2dx,
\end{align*}
we have $s\geq1$, a contradiction. Consequently, we have
 $u_n\rightarrow u$ in $D^{1,2}(\mathbb{R}^3)$.
With the help of the preceding lemmas we can now prove Theorem 1.1.\\
{\textbf{Proof of Theorem 1.1.}}
The theorem follows immediately from Lemma 3.1 and Theorem 2.1.
\section{The case when $\lambda >\lambda_1$}
We assume that $\int_{\mathbb{R}^3}k(x){e_1}^4dx-b(\int_{\mathbb{R}^3}|\nabla u|^2dx)^2 < 0$. Then $e_1\in L^-\bigcap B^-$ and $t(e_1)e_1\in S^+$.
The following lemma plays an important role for establishing the existence of minimizers.
\begin{lem}
Suppose that $\int_{\mathbb{R}^3}\frac{k(x)}{|x|}|e_1|^4dx-b(\int_{\mathbb{R}^3}|\nabla e_1|^2dx)^2 < 0$ . Then there exists $\sigma > 0$ such that $\overline{L^-}\bigcap \overline{B^+}=\emptyset$ whenever $\lambda_1 < \lambda <  \lambda_1+\sigma$.
\end{lem}
\emph{Proof.}
By contradiction, then there exist sequences $\{\lambda_n\}$ and $\{u_n\}$ such that $\|u_n\|=1, \lambda_n\rightarrow \lambda^+_1$ and $$\int_{\mathbb{R}^3}|\nabla u_n|^2-\lambda_n h(x)u_n^2dx \leq 0,$$
$$\int_{\mathbb{R}^3}\frac{k(x)}{|x|}|u_n|^4dx-b(\int_{\mathbb{R}^3}|\nabla u|^2dx)^2\geq 0.$$
Since $u_n$ is bounded, we may assume that $u_n\rightharpoonup u$. We show that $u_n\rightarrow u$ in $D^{1,2}(\mathbb{R}^3)$. Supposing otherwise, then we have $\|u\|<\liminf_{n\rightarrow\infty}\|u_n\|$ and $$\int_{\mathbb{R}^3}|\nabla u|^2-\lambda_1 h(x)|u|^2 < \liminf_{n\rightarrow\infty}\int_{\mathbb{R}^3}|\nabla u_n|^2-\lambda_n h(x)|u_n|^2 \leq 0. $$
This is a contradiction to $\lambda_1$. Hence $u_n\rightarrow u$ in $D^{1,2}(\mathbb{R}^3)$, $\|u\|=1$ and $u=\pm e_1$. To get a contradiction, the cases  $u= e_1$ and $u= -e_1$ are entirely similar, so that we only consider $u= e_1$.
On the other hand, since functional $b(\int_{\mathbb{R}^3}|\nabla u|^2dx)^2$  is continuous on $D^{1,2}(\mathbb{R}^3)$, we have
\begin{equation*}
\lim_{n\rightarrow\infty}b(\int_{\mathbb{R}^3}|\nabla u_n|^2dx)^2 = b(\int_{\mathbb{R}^3}|\nabla e_1|^2dx)^2dx.
\end{equation*}
By Caffarelli-Kohn-Nirenberg inequality (\ref{14}), we have
$$\int_{\mathbb{R}^3}\frac{k(x)}{|x|}|u_n|^4dx \rightarrow \int_{\mathbb{R}^3}\frac{k(x)}{|x|}|e_1|^4dx.$$
Hence,$\int_{\mathbb{R}^3}\frac{k(x)}{|x|}|e_1|^4dx-b(\int_{\mathbb{R}^3}|\nabla e_1|^2dx)^2 \geq 0$ , This is a contradiction to our assumption.
If $\overline{L^-}\cap \overline{B^+}=\emptyset$ is  satisfied, we can get more information on $S$.
\begin{lem}
Suppose that $\overline{L^-}\cap \overline{B^+}=\emptyset$. Then \\
(i)  $S^0 = \{0\}$.\\
(ii) $0\notin \overline{S^-}$ and $S^-$ is closed.\\
(iii) $S^-$  and  $S^+$ are separated. That is $\overline{S^-}\cap \overline{S^+}=\emptyset$.\\
(iv)  $S^+$ is bounded.
\end{lem}
\emph{Proof.}
(i). Suppose $u\in S^0\setminus\{0\}$. Then $\frac{u}{\|u\|}\in L^0 \cap B^0 \subset \overline{L^-}\cap \overline{B^+}=\emptyset$  which is impossible. Hence $S^0 = \{0\}$.\\
(ii). Supposing otherwise, then there exists $\{u_n\}\in S^-$ such that $u_n\rightarrow 0$ in $D^{1,2}(\mathbb{R}^3)$. Hence

\begin{equation}
 0< \int_{\mathbb{R}^3}|\nabla u_n|^2-\lambda h(x)u_n^2dx =
\int_{\mathbb{R}^3}k(x)\frac{|u_n|^4}{|x|}dx-b(\int_{\mathbb{R}^3}|\nabla u_n|^2dx^2)\rightarrow 0.\label{41}
\end{equation}
Let $v_n=\frac{u_n}{\|u_n\|}$. We observe that
\begin{equation}
 0< \int_{\mathbb{R}^3}|\nabla v_n|^2-\lambda h(x)v_n^2dx =
\|u_n\|^2(\int_{\mathbb{R}^3}k(x)\frac{|v_n|^4}{|x|}dx-b(\int_{\mathbb{R}^3}|\nabla v_n|^2dx^2)).\label{42}
\end{equation}
We may assume that $v_n\rightharpoonup v_0$ in $D^{1,2}(\mathbb{R}^3)$. To obtain a contradiction, we divide our proof into three steps:\\
(a)~~~~~~~~~$v_0\neq 0$.\\
(b)~~~~~~~~~$ \frac{v_0}{\|v_0\|}\in \overline{L^-}$.\\
(c)~~~~~~~~~$\frac{v_0}{\|v_0\|}\in \overline{B^+}$.\\
We begin to prove the assertions  (a), (b) and (c).\\
(a). By Caffarelli-Kohn-Nirenberg inequality (\ref{14}), we obtain
 \begin{equation}
 \int_{\mathbb{R}^3}k(x)\frac{|v_n|^4}{|x|}dx \leq C\|v_n\|^4 =C .     \label{43}
\end{equation}
Then we have
 \begin{equation}
 0 <\int_{\mathbb{R}^3}k(x)\frac{|v_n|^4}{|x|}dx-b(\int_{\mathbb{R}^3}|\nabla v_n|^2dx)^2 \leq C.    \label{44}
\end{equation}
Since $u_n\in S^-$ and $\|u_n\| \rightarrow 0$, by (\ref{42}) and (\ref{44}), we have
\begin{equation}
 \lim_{n\rightarrow\infty}\int_{\mathbb{R}^3}|\nabla v_n|^2-\lambda h(x)v_n^2dx =0,   \label{45}
\end{equation}
that is
 $$1=\lim_{n\rightarrow\infty}\int_{\mathbb{R}^3}\lambda h(x)v^2_ndx =\int_{\mathbb{R}^3}\lambda h(x)v^2_0dx.$$
 So $v_0 \neq 0$.\\
(b).By (\ref{45}), we have
 $$\int_{\mathbb{R}^3}|\nabla v_0|^2-\lambda h(x)v_0^2dx\leq\lim_{n\rightarrow\infty}\int_{\mathbb{R}^3}|\nabla v_n|^2-\lambda h(x)v_n^2dx=0.$$
So $ \frac{v_0}{\|v_0\|}\in \overline{L^-}$.\\
(c). According to (\ref{44}), it follows that
\begin{align*}
0\leq\limsup_{n\rightarrow\infty}(\int_{\mathbb{R}^3}k(x)\frac{|v_n|^4}{|x|}dx-b(\int_{\mathbb{R}^3}|\nabla v_n|^2dx)^2) &\leq \limsup_{n\rightarrow\infty}\int_{\mathbb{R}^3}k(x)\frac{|v_n|^4}{|x|}dx -\liminf_{n\rightarrow\infty}b(\int_{\mathbb{R}^3}|\nabla v_n|^2dx)^2 \\
&\leq \limsup_{n\rightarrow\infty}\int_{\mathbb{R}^3}k(x)\frac{|v_n|^4}{|x|}dx-b(\int_{\mathbb{R}^3}|\nabla v_0|^2dx)^2.
\end{align*}
Therefore, according to the definitions of $\alpha_\infty,k_\infty, \nu_0 $, we have
$$\int_{\mathbb{R}^3}k(x)\frac{|v_0|^4}{|x|}dx-b(\int_{\mathbb{R}^3}|\nabla v_0|^2dx)^2\geq -k_\infty\alpha_\infty-k(0)\nu_0\geq0.$$
Hence, $0\notin \overline{S^-}$. Now we are ready to prove that $S^-$  is closed.
 Let $u\in  \overline{S^-}$, then there exist $\{u_n\} \subseteq S^-$ such  that $u_n\rightarrow u$ in $H^1(\mathbb{R}^3)$.
 Since $S$ is closed,  $S^+$ is open, $0\notin \overline{S^-}$  and $S=S^-\cup S^+\cup \{0\}$ , we have $u\in S^-$.  We conclude that $S^-$ is closed.\\\\
(iii). By (i) and (ii),we have
$$\overline{S^-}\cap \overline{S^+}\subseteq S^-\cap(S^+\cup\{0\}=(S^-\cap S^+)\cup(S^-\cap\{0\})=\emptyset.$$
(iv). Suppose that $S^+$ is unbounded. Then there  exists a sequence $\{u_n\}\subset S^+$ such that
\begin{equation*}
\int_{\mathbb{R}^3}|\nabla u_n|^2-\lambda h(x)u_n^2dx =
\int_{\mathbb{R}^3}k(x)\frac{|u_n|^4}{|x|}dx-b(\int_{\mathbb{R}^3}|\nabla u_n|^2dx)^2 \leq 0
\end{equation*}
and $\|u_n\|\rightarrow\infty$.
Let $v_n=\frac{u_n}{\|u_n\|}$. We may assume that $v_n\rightharpoonup v_0$ in $D^{1,2}(\mathbb{R}^3)$.
It is clear that
\begin{equation}
\int_{\mathbb{R}^3}|\nabla v_n|^2-\lambda h(x)v_n^2dx =
\|u_n\|^2(\int_{\mathbb{R}^3}k(x)\frac{|v_n|^4}{|x|}dx-b(\int_{\mathbb{R}^3}|\nabla v_n|^2dx^2))  \label{46}
\end{equation} and also
\begin{equation}
\int_{\mathbb{R}^3}|\nabla v_0|^2-\lambda h(x)v_0^2dx\leq \liminf_{n\rightarrow\infty}\int_{\mathbb{R}^3}|\nabla v_n|^2-\lambda h(x)v_n^2 dx\leq 0. \label{47}
\end{equation}
 We shall adopt the same procedure as  (a),(b),(c) in  (ii).
 Suppose $v_0 =0$. We claim that $v_n\rightarrow 0$ in $D^{1,2}(\mathbb{R}^3)$. Indeed, using (\ref{47}), if $v_n\nrightarrow 0$, we have
 $$ 0 = \int_{\mathbb{R}^3}|\nabla v_0|^2-\lambda h(x)v_0^2dx< \liminf_{n\rightarrow\infty}\int_{\mathbb{R}^3}|\nabla v_n|^2-\lambda h(x)v_n^2 dx\leq 0.$$
 Therefore, $v_n\rightarrow 0$ in $H^1(\mathbb{R}^3)$, which a contradiction to $\|v_n\|=1$. So $\|v_0\|\neq 0$.\\
 By (\ref{47}) and $v_0\neq 0$, we have$$\frac{v_0}{\|v_0\|}\in \overline{L^-}.$$
 On the other hand, according to (\ref{46}), it follows that
 $$ \lim_{n\rightarrow\infty}\int_{\mathbb{R}^3}k(x)\frac{|v_n|^4}{|x|}dx-b(\int_{\mathbb{R}^3}|\nabla v_n|^2dx)^2=0. $$
 According to the definitions of $\alpha_\infty$ and $k_\infty$, then we have
 \begin{align*}
\int_{\mathbb{R}^3}k(x)\frac{|v_0|^4}{|x|}dx-b(\int_{\mathbb{R}^3}|\nabla v_0|^2dx)^2&\geq\limsup_{n\rightarrow\infty}\int_{\mathbb{R}^3}k(x)\frac{|v_n|^4}{|x|}dx -\liminf_{n\rightarrow\infty}b(\int_{\mathbb{R}^3}|\nabla v_n|^2dx)^2-(\alpha_\infty k_\infty+k(0)\nu_0) \\
&\geq\lim_{n\rightarrow\infty}(\int_{\mathbb{R}^3}k(x)\frac{|v_n|^4}{|x|}dx-b(\int_{\mathbb{R}^3}|\nabla v_n|^2dx)^2)-(\alpha_\infty k_\infty+k(0)\nu_0)\\
&=-(\alpha_\infty k_\infty+k(0)\nu_0)\geq0.
\end{align*}
 This means $\frac{v_0}{\|v_0\|}\in \overline{B^+}$. \\
 Summarizing what have proved, (a).$v_0\neq 0$, (b).$ \frac{v_0}{\|v_0\|}\in \overline{L^-}$, (c).$\frac{v_0}{\|v_0\|}\in \overline{B^+}$,
 it is a contradiction to  $\overline{L^-}\cap \overline{B^+}=\emptyset$. Then $S^+$ is bounded.
\begin{lem}
Suppose that $\overline{L^-}\cap \overline{B^+}=\emptyset$. Then \\
(i). Every minimizing sequence of $J(u)$ on $S^-$ is bounded.\\
(ii). $inf_{u\in S^-}J(u)>0$.\\
(iii). There exists a minimizer of $J(u)$ on $S^-$.\\
\end{lem}
\emph{Proof.}
(i). Let $\{u_n\}\in S^-$ be a minimizing sequence for the functional $J(u)$. Then $$J(u_n)=\frac{1}{4}\int_{\mathbb{R}^3}|\nabla u_n|^2-\lambda h(x)u_n^2dx =\frac{1}{4}(\int_{\mathbb{R}^3}k(x)\frac{|u_n|^4}{|x|}dx-b(\int_{\mathbb{R}^3}|\nabla u_n|^2dx)^2)\rightarrow c, $$ where $c>0$.
Similar to the one used  in the proof of (iv) of Lemma 4.2, it is easy to prove that $\{u_n\}$ is bounded in $D^{1,2}(\mathbb{R}^3)$.\\
(ii). It is clear that $J(u)\geq0$ on $S^-$. Suppose $\inf_{w\in S^-}J(w)=0$. Let $\{u_n\}$ be a minimizing sequence. According to (i), it follows  that $\{u_n\}$ is bounded and we may assume $u_n\rightharpoonup u_0$. Clearly, we have
 \begin{equation}
\lim_{n\rightarrow\infty}\int_{\mathbb{R}^3}|\nabla u_n|^2-\lambda h(x)u_n^2dx =0 \label{48}
\end{equation}
and
\begin{equation}
\lim_{n\rightarrow \infty}(\int_{\mathbb{R}^3}k(x)\frac{|u_n|^4}{|x|}dx-b(\int_{\mathbb{R}^3}|\nabla u_n|^2dx)^2) =  0.\label{49}
\end{equation}
Then it is easy to prove that $\frac{u_0}{\|u_0\|}\in \overline{B^+}\cap \overline{L^-}$. Indeed, using the same argument in Lemma 4.2,  by(\ref{48}), it is easy to show that $\frac{u_0}{\|u_0\|}\in \cap \overline{L^-}$ and by (\ref{49}), $\frac{u_0}{\|u_0\|}\in \overline{B^+}$ follows.\\
So  $\text{inf}_{u\in S^-}J(u)>0$.\\
(iii). Let $\{u_n\}$ be a minimizing sequence.
According to (i), it  follows  that $\{u_n\}$ is bounded and we may assume $u_n\rightharpoonup u_0$ in $ H^1(\mathbb{R}^3)$.
Suppose $u_n\nrightarrow u_0$ in $ D^{1,2}(\mathbb{R}^3)$. We get

\begin{eqnarray*}
\int_{\mathbb{R}^3}|\nabla u_0|^2-\lambda h(x)u_0^2dx
&<& \lim_{n\rightarrow\infty}\int_{\mathbb{R}^3}|\nabla u_n|^2-\lambda h(x)u_n^2dx \\
&=&\lim_{n\rightarrow \infty}(\int_{\mathbb{R}^3}k(x)\frac{|u_n|^4}{|x|}dx-b(\int_{\mathbb{R}^3}|\nabla u_n|^2dx)^2) \\
&\leq&\limsup_{n\rightarrow\infty}\int_{\mathbb{R}^3}k(x)\frac{|u_n|^4}{|x|}dx -\liminf_{n\rightarrow\infty}b(\int_{\mathbb{R}^3}|\nabla u_n|^2dx)^2\\
&\leq&\int_{\mathbb{R}^3}k(x)\frac{|u_0|^4}{|x|}dx-b(\int_{\mathbb{R}^3}|\nabla u_0|^2dx)^2+k_\infty\alpha_\infty+k(0)\nu_0\\
&\leq& \int_{\mathbb{R}^3}k(x)\frac{|u_0|^4}{|x|}dx-b(\int_{\mathbb{R}^3}|\nabla u_0|^2dx)^2
\end{eqnarray*}
So there exists a  $ 0<t<1$ such that $tu_0 \in S^- $. Similar to the proof of (ii) in Lemma 3.1, we get a contradiction. \\
Hence $u_n\rightarrow u_0 \neq 0 $ in $ H^1(\mathbb{R}^3)$.
Since $S^-$ is closed, then we have $u_0\in S^-$ and $J(u_0)=\inf_{w\in S^-}J(w).$\\
We are going to the  investigation on $S^+$.
\begin{lem}
Suppose that $\overline{L^-}\cap \overline{B^+}=\emptyset$. Then there exists $0\neq v \in S^+$ such that $J(v) =\inf_{u\in S^+}J(u)$.
\end{lem}
\emph{Proof.}
Due to $\overline{L^-}\cap \overline{B^+}=\emptyset$, $L^-\cap B^-$ as well as $S^+$ must be nonempty. By (iv) of Lemma 4.2, there exists $ M > 0$ such that $\|u\|\leq M$ for all the $u\in S^+$. Using Caffarelli-Kohn-Nirenberg inequality (\ref{14}), it is easy to see that $J(u)$ is bounded from below on $S^+$ and $\inf_{u\in S^+}J(u)< 0 $.
Let $\{u_n\}\in S^+$ be a minimizing sequence for the  functional $J(u)$. Then $$J(u_n)=\frac{1}{4}\int_{\mathbb{R}^3}|\nabla u_n|^2-\lambda h(x)u_n^2dx =\frac{1}{4}(\int_{\mathbb{R}^3}k(x)\frac{|u_n|^4}{|x|}dx-b(\int_{\mathbb{R}^3}|\nabla u_n|^2dx)^2)\rightarrow c, $$ where $c < 0$.
We may assume $u_n\rightharpoonup u_0$ in $ D^{1,2}(\mathbb{R}^3)$. Obviously,
$$\int_{\mathbb{R}^3}|\nabla u_0|^2-\lambda h(x)u_0^2dx \leq \lim_{n\rightarrow\infty}\int_{\mathbb{R}^3}|\nabla u_n|^2-\lambda h(x)u_n^2dx=c<0.$$
So $u_0\neq 0$ and  $\frac{u_0}{\|u_0\|}\in L^-$. Since $\overline{L^-}\cap \overline{B^+}=\emptyset$, then we have $\frac{u_0}{\|u_0\|}\in B^-$. \\
Hence $t(u_0)u_0\in S^+$, where $t(u_0)=[\frac{\int_{\mathbb{R}^3}|\nabla u_0|^2-\lambda h(x)u_0^2dx}{\int_{\mathbb{R}^3}k(x)\frac{|u_0|^4}{|x|}dx-b(\int_{\mathbb{R}^3}|\nabla u_0|^2dx)^2 }]^\frac{1}{2}$.
Suppose that $u_n\nrightarrow u_0$. By an argument similar to the one used in the proof of (iii) in Lemma 4.3, we have
$$\int_{\mathbb{R}^3}|\nabla u_0|^2-\lambda h(x)u_0^2dx <\int_{\mathbb{R}^3}k(x)\frac{|u_0|^4}{|x|}dx-b(\int_{\mathbb{R}^3}|\nabla u_0|^2dx)^2.$$
 Together with $\int_{\mathbb{R}^3}k(x)\frac{|u_0|^4}{|x|}dx-b(\int_{\mathbb{R}^3}|\nabla u_0|^2dx)^2<0$,(since $\frac{u_0}{\|u_0\|}\in B^-$). It implies that $t(u_0)>1$.
This is contrary with
$$J(t(u_0)u_0)<J(u_0)\leq \lim_{n\rightarrow\infty}J(u_n)=\inf_{w\in S^+}J(w).$$
Hence $$u_n \rightarrow u_0$$ in $ D^{1,2}(\mathbb{R}^3)$.
Since the  functional $\Xi(u)=\int_{\mathbb{R}^3}k(x)\frac{|u_n|^4}{|x|}dx-b(\int_{\mathbb{R}^3}|\nabla u_n|^2dx)^2$  is continuous.
Therefore, we have $u_0\in S^+$  and $$J(u_0)=\inf_{w\in S^+}J(u).$$\\
We now turn to the proof of Theorem 1.2.\\
{\textbf{Proof of Theorem1.2}}
According to Lemma 4.1 and the  assumptions of Theorem 1.2, it follows that $\overline{L^-}\cap \overline{B^+} = \emptyset$. We are ready to invoke the
conclusions of Lemma 4.3 and Theorem 2.1.
So, there exists $u_1\in S^-$ which is a critical point of $J(u)$. Clearly, $J(u_1)> 0 $. Employing Lemma 4.4 and Theorem 2.1, there exists  $u_2\in S^+$ which is a
critical point of $J(u)$. Clearly, $J(u_2) < 0 $. We have thus proved   Theorem 1.2.
{\small

\begin{thebibliography}{99}
\bibitem{Alama}
{S. Alama, G. Tarantello}, { \em {On semilinear elliptic equations with indefinite nonlinearities}}, Calc. Var. Partial Differential Equations 1 (1993) 439-475.
\bibitem{Ancona}
{ P. D'Ancona, S. Spagnolo,}{\em{ Global solvability for the degenerate Kirchhoff equation with real analytic data}}, Invent. Math. 108 (2) (1992) 247-262
\bibitem{A-R}
{A. Ambrosetti, P.H. Rabinowitz}, {\em {Dual variational methods in critical point theory and applications}}, J. Funct. Anal. 14 (1973) 349-381.
 \bibitem{Arosio}
{A. Arosio, S. Panizzi,}{\em{ On the well-posedness of the Kirchhoff string}}, Trans. Amer. Math. Soc. 348 (1) (1996) 305-330.
\bibitem{Brown}
{K.~Brown, Y.~Zhang}, {\em {The Nehari manifold for a semilinear elliptic equation with a sign-changing weight function}}, J. Differential Equations 193 (2003) no. 2, 481-49.
\bibitem{C-K-N}
{L. Caffarelli, R. Kohn and L. Nirenberg,}{{\em First order interpolation inequalities with weights}}, Compos. Math. 53 (1984), 259-275.
\bibitem{C-V}
 {G.~Cerami, G~Vaira}, {\em {Positive solution for some non-autonomous Schr\"{o}dinger-Poisson systems}}, J. Differential Equations 248 (2010) 521-543.
 \bibitem{Chabrowski-Costa}
{J. Chabrowski,  D.G. Costa}, { \em  {On a class of Schr\"{o}dinger-type equations
with indefinite weight functions}},Communications in Partial Differential Equations, 33 1368-1394, 2008.
\bibitem{Chabrowski1}
{J.~Chabrowski},{\em{ Weak Convergence Methods for Semilinear Elliptic
Equations. Singapore}}, World Scientific 1999.
\bibitem{Chabrowski2}
{J. Chabrowski, D.G. Costa}, {\em {On existence of positive solutions for a class of Caffarelli-Kohn-Nirenberg type equations}},
 Colloq. Math. 120 (2010) 43-62.
\bibitem{C-K-W}
{C.Chen, Y. Kuo, T. Wu},{\em {The Nehari manifold for a Kirchhoff type problem involving sign-changing weight functions}}, J. Differential Equations 250 (2011), no. 4, 1876-1908.
\bibitem{Chen}
{J.~Chen}, {{\em Multiple positive solutions to a class of Kirchhoff equation on R3 with indefinite nonlinearity}}, Nonlinear Anal. 96 (2014), 134-145.
%\bibitem{Costa-T}
% {D.G. Costa, H. Tehrani},{ \em { Existence of positive solutions for a class of indefinite elliptic problems in $R^N$}},
%Calc. Var. Partial Differential Equations 13 (2001) 159-189.
%\bibitem{Mugnai1}
%{T. D'Aprile, D. Mugnai}, {\em{ Non-existence results for the coupled Klein-Gordon-Maxwell equations}}, Adv. Nonlinear Stud.
%4 (2004) 307-322.
%\bibitem{Chen1}
%{J. Chen}
%{\em{Multiple positive solutions for a class of nonlinear elliptic equations}},
%J. Math. Anal. Appl. 295 (2004), no. 2, 341-354.
\bibitem{Chen-li}
{S. Chen,  L. Li} {\em {Multiple solutions for the nonhomogeneous Kirchhoff equation on $\mathbb{R}^3$}}.
 Nonlinear Anal. Real World Appl. 14 (2013), no. 3, 1477-1486.
\bibitem{He}
{X. M. He, W. M. Zou},{\em{ Existence and concentration behavior of positive solutions for a Kirchhoff equation in $R^3$}}, J. Differential Equations 252 (2012) 1813-1834.
\bibitem{Li-Ye}
{G. Li, H. Ye}, {{\em Existence of positive ground state solutions for the nonlinear Kirchhoff type equations in $R^3$}}, J. Differential Equations 257 (2014), no. 2, 566-600
\bibitem{lions1}
{J. L. Lions}, {\em {On some questions in boundary value problems of mathematical physics, in: Contemporary Development in Continuum Mechanics and Partial Differential Equations}} in: North-Holland Math. Stud., vol. 30, North-Holland, Amsterdam, New York, 1978,  284-346.
\bibitem{Liu-Sun}
{X. Liu; Y. Sun,} {\em{ Multiple positive solutions for Kirchhoff type problems with singularity}},
  Commun. Pure Appl. Anal. 12 (2013), no. 2, 721-733.
%  \bibitem{Smets}
%{D. Smets,} {\em {Nonlinear Schr\"{o}dinger equations with Hardy potential and critical nonlinearities,}}
%Trans. Amer. Math. Soc. 357 (2005), no. 7, 2909-2938.

\bibitem{Willem}
{M. Willem}, {\em{Minimax Theorems}}, Progr. Nonlinear Differential Equations Appl., vol. 24, Birkh\"{a}user Boston, Inc., Boston, MA, 1996.

 \bibitem{Wang}
{J. Wang, L.Tian, J. Xu, F. Zhang}, {\em{Multiplicity and concentration of positive solutions for a Kirchhoff type problem with critical growth}}, J. Differential Equations 253 (2012) 2314-2351.
\bibitem{Wu}
{X. Wu,} {{\em Existence of nontrivial solutions and high energy solutions for Schr\"{o}dinger-Kirchhoff-type equations in $R^N$}}, Nonlinear Anal. Real World Appl. 12 (2011) 1278-1287.

\end{thebibliography}
\end{document}